\documentclass[a4paper,12pt]{article}
\usepackage{amsmath}
\usepackage{color}
\usepackage{amssymb}
\usepackage{latexsym}
\usepackage{amsthm}
\usepackage{url}
\usepackage{multirow}
\usepackage{booktabs}
\usepackage{enumitem}
\usepackage{verbatim}
\numberwithin{equation}{section}

\def\d{\displaystyle}

\def\e{\varepsilon}

\def\R{{\bf R}}
\def\N{{\bf N}}

\def\q{\quad}

\newtheorem{thm}{Theorem}[section]

\newtheorem{rem}{Remark}[section]

\title{Slicing method for nonlinear integral inequalities\\
related to critical nonlinear wave equations}
\author{
Takiko Sasaki
\footnote{
Department of Mathematical Engineering, Faculty of Engineering, Musashino University,
3-3-3 Ariake, Koto-ku, Tokyo 135-8181, Japan.
e-mail: t-sasaki@musashino-u.ac.jp.
}
, Kerun Shao
\footnote{School of Mathematical Sciences, Zhejiang University, Hangzhou 310058, P. R. China,\newline
email: shaokr@163.com.}
, Hiroyuki Takamura
\footnote{Mathematical Institute,
Tohoku University,
Aoba, Sendai 980-8578, Japan.\newline
e-mail: hiroyuki.takamura.a1@tohoku.ac.jp.}
}
\date{
\[
\begin{array}{ll}
\mbox{\footnotesize{\bf Keywords:}}
& \mbox{\footnotesize nonlinear wave equation, blow-up, critical power,}\\
& \mbox{\footnotesize ordinary differential inequality, lifespan, slicing method}\\
\mbox{\footnotesize{\bf MSC2020:}}
& \mbox{\footnotesize primary 35L71, secondary 35B44}\\
\end{array}
\]
}

\begin{document}
\maketitle

\begin{abstract}
This paper is devoted to a simple and short proof on the sharp upper bound of lifespan
of classical solutions to wave equations with the critical power nonlinearities of spatial derivatives of the unknown function.
Such a proof is so-called \lq\lq slicing method",
which may help us to extend the result for various equations and systems.
\end{abstract}


\section{Introduction}
In this paper, we are focusing on the following system of integral inequalities
for unknown functions $H\in C([R,T))$: 
\begin{equation}
\label{NII}
\left\{
\begin{array}{ll}
\d H(t)\geq At^a(\log t)^{-b}\left(\log\frac{t}{R}\right)^c,
&\q t\in [R, T),\\
\d H(t)\geq B(\log t)^x\int_{R}^{t}ds\int_{R}^{s}r^{y}\left(\log \frac{r}{R}\right)^z|H(r)|^p dr,
&\q t\in [R, T),
\end{array} 
\right.
\end{equation}
where all $A>0$, $B>0$ and $T>R>1$ are constants.
We assume that exponents $a,b,c,x,y,z,p$ satisfy
\begin{equation}
\label{exponents}
\left\{
\begin{array}{l}
\d p>1,\ a\leq 1 ,\ b\geq \max\left\{0,\frac{x}{p-1}\right\},  \\
\d y+pa=-1, \ z+cp>-1, \ z+cp\geq c-1.
\end{array}
\right.
\end{equation}
When
\[
a=1,\ b=0,\ c=1,\ x=-p,\ y=-p-1,\ z=1,
\]
\eqref{NII} can be found in Shao, Takamura and Wang \cite[Lemma 3.1]{STW2025},
which leads to blow-up, as well as the optimal upper bound of the lifespan,
of classical solutions to wave equations with critical power nonlinearities of spatial derivative-type;
\begin{equation}
\label{NLW}
\left\{
\begin{array}{ll}
u_{tt}-\Delta_xu=|\nabla_xu|^{(n+1)/(n-1)} & \mbox{in}\ \R^n\times(0,T),\\
u(x,0)=\e f(x),\ u_t(x,0)=\e g(x), & x\in\R^n,
\end{array}
\right.
\end{equation}
where $n\ge2$, $f,g\in C_0^\infty(\R^n)$ and $0<\e\ll1$.
See Introduction of \cite{STW2025} and references therein for the background and criticality
of the problem \eqref{NLW}.
The lifespan $T(\e)$ is defined as the maximal existence time of the solution and \cite{STW2025} shows
\begin{equation}
\label{upper}
T(\e)\le\exp\left(C \e^{-2/(n-1)}\right),
\end{equation}
where $C$ is a positive constant independent of $\e$.
We remark that there is no critical case in one space dimension; see Sasaki, Takamatsu and Takamura \cite{STT2023} on the equation
$u_{tt}-u_{xx}=|u_x|^p\ (p>1)$,
and see Haruyama and Takamura \cite{HT2025} for its application to quasilinear versions.
If $\nabla_x u$ in the nonlinear term in \eqref{NLW} is replaced with $u_t$,
the proof of the blow-up is completely different,
and is reduced to a comparison with one first-order ordinary differential equation; see Zhou \cite{Zhou2001} for its details which is available in all space dimensions.

\par
The purpose of this paper is not only to give a short and simple proof of Lemma 3.1 in \cite{STW2025},
which is the key lemma to show \eqref{upper}, but also to extend the results under more general setting on the exponents as in \eqref{exponents} by making use of the so-called \lq\lq slicing method",
which enable us to prove blow-up of $H$ in \eqref{NII}
by iteration argument with logarithmic terms.
We note that the proof of Lemma 3.1 in \cite{STW2025} has 5 pages with relatively complicated arguments.
But one can see below that we have 2 pages with simple computations even for more general cases.
The iteration argument was first introduced to the point-wise estimate of the solution
to semilinear wave equations $u_{tt}-\Delta_x u=|u|^p\ (p>1)$ in three space dimensions
by John \cite{John1979} to show blow-up in the sub-critical cases.
And it was first combined with the functional method by Lai and Takamura \cite{LT2018}
to show the sharp upper bound of the lifespan in the sub-critical case of energy solutions to
semilinear damped wave equations in the scattering case $u_{tt}-\Delta_x u+\mu u_t/(1+t)^\beta=|u|^p\ 
(p>1,\mu>0,\beta>1)$ in all space dimensions greater than 1.
The critical cases cannot be covered by these methods due to the speciality of this nonlinear term $|u|^p$.
See Ikeda, Sobajima and Wakasa \cite{ISW2019} including references therein
for the history and the latest proof of the critical case of semilinear wave equations,
and see Wakasa and Yordanov \cite{WY2019damped} for those of damped wave equations.
Roughly speaking, as shown in the table below, we can summarize critical exponents for nonlinear wave equations with small data,
including some open parts though.
\begin{center}
\begin{tabular}{|c|c|}
\hline
equation  & critical exponent \\
\hline
\hline
$u_{tt}-\Delta_x u=|u_t|^p$ & $p_G(n)$\\
\hline
$u_{tt}-\Delta_x u=|\nabla_x u|^p$ & $p_G(n)$ \\
\hline
$u_{tt}-\Delta_x u=|u|^p$ & $p_S(n)$\\
\hline
\end{tabular}
\end{center}
Here, so-called Glassey exponent $p_G(n)$ and Strauss exponent $p_S(n)$ are defined by
\[
p_G(n):=\frac{n+1}{n-1},\quad p_S(n):=\frac{n+1+\sqrt{n^2+10n-7}}{2(n-1)}\quad(n\ge2).
\]
Since the blow-up results for critical and sub-critical exponents are established, the estimate of the lifespan $T(\e)$ becomes our primary concern.
\par
Finally, we point out that the slicing method was first introduced by Agemi, Kurokawa and Takamura \cite{AKT2000}
to obtain the optimal upper bound of the lifespan of the solutions to weakly coupled systems
of semilinear wave equations in three space dimensions.
Later, this method has been applied to various equations.
For example, see Wakasa and Yordanov \cite{WY2019} for variable coefficient case,
and  Kitamura, Takamura and Wakasa \cite{KTW2023} for weighted nonlinearities in one dimension.
Therefore, one can expect to extend our result to more various nonlinearities than \eqref{NLW},
and also to weakly coupled systems easily.

\par
This paper is organized as the theorem and its proof are in the next section.
The concluding remarks are added at the end of this paper.
This work was carried out when the first author had employed in a cross-appointment system
between Tohoku Univ. and Musashino Univ. on Apr. 2020 - Mar. 2025,
and during the second author' s stay in Tohoku Univ. on Nov. 2024 - Apr. 2025.



\section{Theorem and proof}

We shall prove the following theorem.

\begin{thm}
\label{thm:main}
Let $H\in C([R,T))$ be a solution to \eqref{NII}
Under assumption \eqref{exponents}, $T$ has to satisfies
\begin{equation}
\label{lifespan}
T\leq\exp\left(\max\left\{2\log R_\infty,(A^{-1}D)^{\frac{p-1}{x+z+1+(c-b)(p-1)}}\right\}\right),
\end{equation}
where
\[
\begin{array}{l}
\d R_\infty=R\prod_{k=1}^{\infty}(1+2^{-k}),\\
\d D=2^{\left(c+\frac{p+(z+1)(p-1)}{(p-1)^2}\right)}p^{\frac{p}{(p-1)^2}}\left(\frac{1}{B}\max\left\{c+\frac{z+1}{p-1},
c+\frac{z+1}{p}\right\}\right)^{\frac{1}{p-1}}>0.
\end{array}
\]
\end{thm}

\begin{rem}
When $a=1$, $b=0$, $c=1$, $x=-p$, $y=-p-1$, $z=1$,
the estimate \eqref{lifespan} becomes
$$
T\leq\exp\left(\max\left\{2\log R_\infty,(A^{-1}D)^{p-1}\right\}\right)
$$
and yields the result in Lemma 3.1 of \cite{STW2025}.
Then, using modified Rammaha's functionals, we can obtain estimates \eqref{upper} for the solutions to equation \eqref{NLW}; see \cite{STW2025} for details.
\end{rem}

\par\noindent
\begin{proof}[\bf\itshape Proof of Theorem \ref{thm:main}]
Substituting $H$ in the right hand side of the first line in \eqref{NII} into the right hand side of the second line in \eqref{NII},
we have that for $t\ge R$
\[
\begin{array}{ll}
H(t)
&\d\geq A^pB(\log t)^x\int_{R}^{t}ds\int_{R}^{s}r^{y+pa}(\log r)^{-pb}\left(\log \frac{r}{R}\right)^{pc+z}dr\\
&\d\geq A^pB(\log t)^{-(pb-x)}\int_{R}^{t}ds\int_{R}^{s}r^{-1}\left(\log \frac{r}{R}\right)^{pc+z}dr\\
&\d=\frac{A^pB}{pc+z+1}(\log t)^{-(pb-x)}\int_{R}^{t}\left(\log\frac{s}{R}\right)^{pc+z+1}ds
\end{array}
\]
because $b\ge0$, $y+pa=-1$ and $pc+z+1>0$ by assumption \eqref{exponents}.
Then, for $\delta>0$ and $t\ge(1+\delta)R$, $H(t)$ can be estimated from below as 
\begin{equation}
\label{eq-iteration-inequal}
\begin{aligned}
H(t)&\geq \frac{A^pB}{pc+z+1}(\log t)^{-(pb-x)}\int_{\frac{t}{1+\delta}}^{t}\left(\log\frac{s}{R}\right)^{pc+z+1}ds\\
&\geq \frac{\delta A^pB}{(1+\delta)(pc+z+1)}t(\log t)^{-(pb-x)}\left(\log\frac{t}{(1+\delta)R}\right)^{pc+z+1}\\
&\geq \frac{\delta A^pB}{(1+\delta)(pc+z+1)}t^a(\log t)^{-(pb-x)}\left(\log\frac{t}{(1+\delta)R}\right)^{pc+z+1}
\end{aligned}
\end{equation}
since $a\leq1$.

\par
Now, set $\delta_j:=2^{-j}$ for $j\in \N$
and define sequences $\{b_j\}$, $\{c_j\}$, $\{R_j\}$, and $\{A_j\}$ by
\[
\left\{
\begin{array}{ll}
b_{j+1}= pb_j-x, & b_0=b,\\
c_{j+1}=pc_j+z+1,& c_0=c,\\
R_{j+1}=(1+\delta_{j})R_{j},& R_0=R,\\
A_{j+1}=\d\frac{\delta_{j} A_j^pB}{(1+\delta_{j})(pc_j+z+1)},& A_0=A.
\end{array}
\right.
\]
Then, by direct calculations, it follows that for $j\in \N$ 
\[
\begin{array}{l}
\d b_{j}=p^{j}\left(b-\frac{x}{p-1}\right)+\frac{x}{p-1},
\quad c_{j}=p^{j}\left(c+\frac{z+1}{p-1}\right)-\frac{z+1}{p-1},\\
\d R_{j+1}=R\prod_{k=0}^{j}(1+2^{-k}).
\end{array}
\]
Since
\[
b\geq \max\left\{0,\frac{x}{p-1}\right\},\ z+cp>-1,\ \mbox{and}\ z+cp\geq c-1,
\]
we have $b_j\geq0$ and $c_{j+1}>0$ for all $j\in \N$. 
Thus, according to the calculation in \eqref{eq-iteration-inequal}, we deduce that
\[
H(t)\geq A_jt^a(\log t)^{-b_j}\left(\log\frac{t}{R_j}\right)^{c_j}\quad\mbox{for}\ t\geq R_j
\]
implies
\[
H(t)\geq A_{j+1}t^a(\log t)^{-b_{j+1}}\left(\log\frac{t}{R_{j+1}}\right)^{c_{j+1}}\quad\mbox{for}\ t\geq R_{j+1}.
\]
Note that
\begin{equation}
\label{eq-Aj-iteration-esitmate}
\begin{aligned}
\log A_{j+1}&=p\log A_{j}+\log B-\log(1+2^{j})-\log c_{j+1}\\
&\geq p\log A_{j}-(j+1)\log (2p)-\log C,
\end{aligned}
\end{equation}
where
\[
C:=\frac{1}{B}\max\left\{c+\frac{z+1}{p-1},c+\frac{z+1}{p}\right\}>0.
\]
Then, inequality \eqref{eq-Aj-iteration-esitmate} is equivalent to
\[
\begin{aligned}
&\log A_{j+1}-(j+1)\frac{\log (2p)}{p-1}-\frac{\log\left((2p)^{\frac{p}{p-1}}C\right)}{p-1}\\
&\geq p\left(\log A_{j}-j\frac{\log (2p)}{p-1}-\frac{\log\left((2p)^{\frac{p}{p-1}}C\right)}{p-1}\right).
\end{aligned}
\]
Hence, we obtain that
\[
\log A_{j}\geq p^{j}\log\left(\frac{A}{(2p)^{\frac{p}{(p-1)^2}}C^{\frac{1}{p-1}}}\right)+j\frac{\log (2p)}{p-1}+\log\left((2p)^{\frac{p}{(p-1)^2}}C^{\frac{1}{p-1}}\right)
\]
for $j\in \N$.

Let $R_\infty$ denote $\d R\prod_{k=1}^{\infty}(1+2^{-k})$. 
It follows that
\[
\begin{aligned}
H(t)&\geq A_jt^a(\log t)^{-b_j}\left(\log\frac{t}{R_j}\right)^{c_j}\\
&\geq (2p)^{\frac{p+(p-1)j}{(p-1)^2}}C^{\frac{1}{p-1}}\left(\frac{A}{(2p)^{\frac{p}{(p-1)^2}}C^{\frac{1}{p-1}}}\right)^{p^j}t^a(\log t)^{-b_j}\left(\log\frac{t}{R_\infty}\right)^{c_j}
\end{aligned}
\]
for $t\geq R_\infty$ and $j\in \N$.
Then, for $t\geq R_\infty^2$ and $j\in \N$, we have
\begin{equation}
\label{eq-inquality-geqtildeR2}
\begin{aligned}
H(t)&\geq (2p)^{\frac{p+(p-1)j}{(p-1)^2}}C^{\frac{1}{p-1}}\left(\frac{A}{(2p)^{\frac{p}{(p-1)^2}}C^{\frac{1}{p-1}}}\right)^{p^j}t^a(\log t)^{-b_j}\left(\frac{1}{2}\log t\right)^{c_j}\\
&=2^{\frac{z+1}{p-1}}(2p)^{\frac{p+(p-1)j}{(p-1)^2}}C^{\frac{1}{p-1}}t^a\left(\log t\right)^{-\frac{x+z+1}{p-1}}\left(\frac{A(\log t)^{\left(c-b+\frac{x+z+1}{p-1}\right)}}{2^{\left(c+\frac{z+1}{p-1}\right)}(2p)^{\frac{p}{(p-1)^2}}C^{\frac{1}{p-1}}}\right)^{p^j}.
\end{aligned}
\end{equation}

\par
If $T$ would satisfy
\[
\frac{A(\log T)^{\left(c-b+\frac{x+z+1}{p-1}\right)}}{2^{\left(c+\frac{z+1}{p-1}\right)}(2p)^{\frac{p}{(p-1)^2}}C^{\frac{1}{p-1}}}>1,
\]
$H(T)$ cannot be finite as $j\rightarrow\infty$.
Therefore, $T$ should satisfy
\[
T\leq\exp\left(\max\left\{2\log R_\infty,(A^{-1}D)^{\frac{p-1}{x+z+1+(c-b)(p-1)}}\right\}\right),
\]
where
\[
D:=2^{\left(c+\frac{p+(z+1)(p-1)}{(p-1)^2}\right)}p^{\frac{p}{(p-1)^2}}C^{\frac{1}{p-1}}>0.\qedhere
\]	
\end{proof}


\section{Concluding remarks}
\par
We have shown a simple proof to obtain the sharp upper bound of the lifespan of the solution to \eqref{NLW}.
\lq\lq Slicing method" may help us to have an application to weakly coupled systems such as
\[
\left\{
\begin{array}{l}
u_{tt}-\Delta u=|\nabla_xv|^p,\\
v_{tt}-\Delta v=|\nabla_xu|^q
\end{array}
\right.
\quad(p,q>1)
\]
with a critical exponents $(p,q)$ as well as its damped version.
For references on the nonlinearities in which $\nabla_xu$ is replaced with $u_t$,
see Ikeda, Sobajima and Wakasa \cite{ISW2019}.
Also it is interesting to see \lq\lq combined effect" by new model,
\[
u_{tt}-\Delta u=|\nabla_xu|^p+|u|^q\quad(p,q>1).
\]
See Introduction of Kido, Sasaki, Takamatsu and Takamura \cite{KSTT2024} for all the references
on the combined effect for nonlinear wave equations. 

\par
In a recent hot topic in this research area, one may have an application of our work to
\lq\lq modulus continuity" such as
\[
u_{tt}-\Delta u=|\nabla_xu|^{p_G(n)}\mu(|\nabla_xu|),
\]
where $\mu$ is some function which is weaker than any power but contributes the total
integrability of a function $|s|^{p_G(n)}\mu(s)$.
For references on the modulus continuity, 
see Chen and Reissig \cite{CR2024} or Wang and Zhang \cite{WZ}
on which $\nabla_xu$ and $p_G(n)$ are replaced with $u$ and $p_S(n)$,
and see Chen and Palmieri \cite{CP} or Shao \cite{S}
on which $\nabla_xu$ is replaced with $u_t$.
In this way, we have many possibilities to apply our result.

\section*{Acknowledgement}
\par
The first and third authors are partially supported
by the Grant-in-Aid for Scientific Research (A) (No.22H00097) and (C) (No.24K06819), 
Japan Society for the Promotion of Science.
The second author is partially supported
by China Scholarship Council (No. 202406320284).


\bibliographystyle{plain}

\end{document}